\documentclass[11pt,a4paper]{article}

\usepackage{amsmath,amssymb}
\newtheorem{theorem}{Theorem}[section]

\newtheorem{definition}{Definition}[section]
\newtheorem{remark}{Remark}[section]
\numberwithin{equation}{section}
\newtheorem{example}{Example}[section]
\newenvironment{proof}[1][Proof]{\noindent\textbf{#1.} }{\hfill {$\Box$}}
\numberwithin{equation}{section}

\begin{document}

\title{\textsc{Nonuniform Behaviors for Skew-Evolution Semiflows in Banach Spaces }}
\author{\textsc{Codru\c{t}a Stoica} \ \ \textsc{Mihail Megan}}
\date{}
\maketitle

{\footnotesize \noindent \textbf{Abstract.} The paper emphasizes
some asymptotic behaviors for skew-evolution semiflows in Banach
spaces. These are defined by means of evolution semiflows and
evolution cocycles. Some characterizations which generalize
classical results are also provided. The approach is from
nonuniform point of view. }

{\footnotesize \vspace{3mm} }


{\footnotesize \vspace{2mm} }


\section{Preliminaries}

The study of asymptotic properties, such as exponential dichotomy
and exponential trichotomy, considered basic concepts that appear
in the theory of dynamical systems, plays an important role in the
study of stable, instable and central manifolds. Some of the
original results concerning stability and instability were
published in \cite{MeStBu_CJM}, \cite{St_arXiv_1} and
\cite{St_arXiv_3} for a particular case of skew-evolution
semiflows defined by means of semiflows and cocycles. In this very
case was also defined and characterized the trichotomy on Banach
spaces in \cite{MeSt_CAMA}.

Concerning previous results, C. Bu\c{s}e presents in
\cite{Bu_RSMUPT} the nonuniform exponential stability for
evolutionary processes. Caracterizations for the nonuniform
exponential instability for evolution operators on Banach spaces
were obtained by M. Megan, A.L. Sasu and B. Sasu in
\cite{MeSaSa_GM}. The study of the nonuniform exponential
dichotomy for evolution families was emphasized by P. Preda and M.
Megan in \cite{PrMe_BAMS} and for evolution operators, also in
Banach spaces, by M. Megan, A.L. Sasu and B. Sasu in
\cite{MeSaSa_IEOT02}. Other asymptotic properties for evolution
families were studied by the same authors in the nonuniform
setting in \cite{MeSaSa_MIR}.

In this paper we extend the asymptotic
properties of exponential dichotomy and trichotomy for the newly
introduced concept of skew-evolution semiflows defined on Banach spaces, which can be
considered generalizations for evolution operators and skew-product semiflows. The results concerning the nonuniform exponential trichotomy are
generalizations of some Theorems proved for evolution operators in \cite{MeSt_IEOT}.

\section{Notations. Definitions. Examples}

Let us consider a metric space $(X,d)$, a Banach space $V$ and
$\mathcal{B}(V)$ the space of all bounded linear operators from
$V$ into itself. Let $V^{*}$ be the topological dual of $V$. We denote the sets $T =\left\{(t,t_{0})\in
\mathbb{R}^{2}, \ t\geq t_{0}\geq 0\right\}$ and $Y=X\times V$.
Let $P:Y\rightarrow Y$ be a projector given by $P(x,v)=(x,P(x)v)$,
where $P(x)$ is a projection on $Y_{x}=\{x\}\times V$, $x\in X$.

\begin{definition}\rm\label{def_sfl_ev}
A mapping $\varphi: T\times  X\rightarrow  X$ is called
\emph{evolution semiflow} on $ X$ if following relations hold:

$(s_{1})$ $\varphi(t,t,x)=x, \ \forall (t,x)\in
\mathbb{R}_{+}\times X$

$(s_{2})$ $\varphi(t,s,\varphi(s,t_{0},x))=\varphi(t,t_{0},x),
\forall (t,s),(s,t_{0})\in T, x\in X$.
\end{definition}

\begin{definition}\rm\label{def_aplcoc_ev}
A mapping $\Phi: T\times  X\rightarrow \mathcal{B}(V)$ is called
\emph{evolution cocycle} over an evolution semiflow $\varphi$ if:

$(c_{1})$ $\Phi(t,t,x)=I$, the identity operator on $V$, $\forall
(t,x)\in \mathbb{R}_{+}\times X$

$(c_{2})$
$\Phi(t,s,\varphi(s,t_{0},x))\Phi(s,t_{0},x)=\Phi(t,t_{0},x),\forall
(t,s),(s,t_{0})\in T, x\in X$.
\end{definition}

\begin{definition}\rm\label{def_coc_ev_1}
The mapping $C: T\times Y\rightarrow Y$ defined by the relation
$C(t,s,x,v)=(\varphi(t,s,x),\Phi(t,s,x)v)$, where $\Phi$ is an
evolution cocycle over an evolution semiflow $\varphi$, is called
\emph{skew-evolution semiflow} on $Y$.
\end{definition}

\begin{example}\rm\label{ex_ses}
We denote $\mathcal{C}=\mathcal{C}(\mathbb{R}_{+},\mathbb{R}_{+})$
the set of all continuous functions $x:\mathbb{R}_{+}\rightarrow
\mathbb{R}_{+}$, endowed with the topology of uniform convergence on compact subsets of $\mathbb{R}_{+}$
and which is metrizable by means of the distance
\begin{equation*}
d(x,y)=\sum_{n=1}^{\infty}\frac{1}{2^{n}}\frac{d_{n}(x,y)}{1+d_{n}(x,y)},
\ \textrm{where} \ d_{n}(x,y)=\underset{t\in
[0,n]}\sup{|x(t)-y(t)|}.
\end{equation*}
If $x\in \mathcal{C}$ then for all $t\in
\mathbb{R}_{+}$ we denote $x_{t}(s)=x(t+s)$, $x_{t}\in
\mathcal{C}$. Let $ X$ be the closure in $\mathcal{C}$ of the set $\{f_{t},t\in \mathbb{R}_{+}\}$, where
$f:\mathbb{R}_{+}\rightarrow \mathbb{R}_{+}^{*}$ is a decreasing function with the property
$\underset{t\rightarrow \infty}\lim{f(t)}=l>0$. Then $( X,d)$
is a metric space and the mapping
\[
\varphi: T\times X\rightarrow  X, \
\varphi(t,s,x)=x_{t-s}
\]
is an evolution semiflow on $X$.

We consider the Banach space $V=\mathbb{R}^{n}$, $n\geq 1$, with the
norm
\[
\left\Vert
(v_{1},...,v_{n})\right\Vert=|v_{1}|+...+|v_{n}|.
\]
The mapping $\Phi: T\times  X\rightarrow \mathcal{B}(V)$ given by
\[
\Phi(t,s,x)(v_{1},...,v_{n})=\left(
e^{\alpha_{1}\int_{s}^{t}x(\tau-s)d\tau}v_{1},...,e^{\alpha_{n}\int_{s}^{t}x(\tau-s)d\tau}v_{n}\right),
\]
where $\alpha=(\alpha_{1},...,\alpha_{n})\in\mathbb{R}^{n}$, is an evolution cocycle over $\varphi$ and $C=(\varphi,\Phi)$ is a
skew-evolution semiflow on $Y$.
\end{example}

An interesting class of skew-evolution semiflows, useful to describe asymptotic properties, is given by

\begin{example}\rm\label{cevshift}
Let us consider a skew-evolution semiflow $C=(\varphi, \Phi)$
and a parameter $\lambda \in \mathbb{R}$. We define the application
\begin{equation}\label{relcevshift}
\Phi_{\lambda}: T\times  X\rightarrow \mathcal{B}(V), \
\Phi_{\lambda}(t,t_{0},x)=e^{-\lambda(t-t_{0})}\Phi(t,t_{0},x).
\end{equation}
It is to remark that $C_{\lambda}=(\varphi, \Phi_{\lambda})$, where $\Phi_{\lambda}$
verifies the conditions of Definition \ref{def_aplcoc_ev}, is a skew-evolution semiflow and it
will be defined as the \emph{$\lambda$-shift} skew-evolution semiflow on $Y$.
\end{example}

\begin{definition}\rm\label{tc_tm}\label{def_startaremas}
A skew-evolution semiflow $C =(\varphi,\Phi)$ is said to be

$(sm)$ \emph{strongly measurable} if for all $(t_{0},x,v)\in T\times Y$ the mapping
$s\mapsto\left\Vert\Phi(s,t_{0},x)v\right\Vert$ is measurable on $[t_{0},\infty)$.

$(ssm)$ \emph{$*$-strongly measurable} if for all
$(t,t_{0},x,v^{*})\in T\times X\times V^{*}$ the mapping
$s\mapsto\left\Vert\Phi(t,s,\varphi(s,t_{0},x))^{*}v^{*}\right\Vert$ is measurable on $[t_{0},t]$.
\end{definition}

\begin{definition}\rm\label{def_neg}
The skew-evolution semiflow $C$ is said to have \emph{exponential
growth} if there exist some applications $M$,
$\omega:\mathbb{R}_{+}\rightarrow\mathbb{R}_{+}^{*}$ such that
\begin{equation}
\left\Vert \Phi(t,t_{0},x)v\right\Vert \leq M(s)e^{\omega(s)
(t-s)}\left\Vert \Phi(s,t_{0},x)v\right\Vert,
\end{equation}%
for all $(t,s),(s,t_{0})\in  T$ and all $(x,v)\in Y$.
\end{definition}

\begin{remark}\rm\label{obs_shift}
Let $C=(\varphi,\Phi)$ be a skew-evolution semiflow with
exponential growth and let $C_{-\alpha}=(\varphi,\Phi_{-\alpha})$,
$\alpha>0$, be the \emph{$-\alpha$-shift} skew-evolution semiflow,
where the evolution cocycle $\Phi_{-\alpha}$ is given by relation
(\ref{relcevshift}). We have
\[
\left\Vert\Phi_{-\alpha}(t,t_{0},x)v\right\Vert=e^{\alpha(t-t_{0})}
\left\Vert\Phi(t,t_{0},x)v\right\Vert \leq
M(t_{0})e^{[\alpha+\omega(t_{0})](t-t_{0})}\left\Vert v\right\Vert
\]
for all $(t_{0},x,v)\in\mathbb{R}_{+}\times Y$, where the
functions $M$ and $\omega$ are given by Definition \ref{def_neg}.
Hence, $C_{-\alpha}$ has also exponential growth.
\end{remark}

\begin{definition}\rm\label{def_nedc}
The skew-evolution semiflow $C$ is said to have \emph{exponential
decay} if there exist some applications $M$,
$\omega:\mathbb{R}_{+}\rightarrow\mathbb{R}_{+}^{*}$ such that
\begin{equation}
\left\Vert \Phi(s,t_{0},x)v\right\Vert \leq M(t)e^{\omega(t)
(t-s)}\left\Vert \Phi(t,t_{0},x)v\right\Vert,
\end{equation}%
for all $(t,s),(s,t_{0})\in  T$ and all $(x,v)\in Y$.
\end{definition}

\begin{remark}\rm
Let $C=(\varphi,\Phi)$ be a skew-evolution semiflow with exponential decay
and let
$C_{\alpha}=(\varphi,\Phi_{\alpha})$, $\alpha>0$, be the $\alpha$-shifted skew-evolution semiflow, where the evolution cocycle
$\Phi_{\alpha}$ is given by relation
(\ref{relcevshift}). Following relations
\[
\left\Vert\Phi_{\alpha}(s,t_{0},x)v\right\Vert=e^{-\alpha(s-t_{0})}
\left\Vert\Phi(s,t_{0},x)v\right\Vert \leq
M(t)e^{[\omega(t)+\alpha](t-s)}\left\Vert
\Phi_{\alpha}(t,t_{0},x)v\right\Vert
\]
hold for all $(t,s),(s,t_{0})\in T$ and all $(x,v)\in Y$, where
the functions $M$ and $\omega$ are given by Definition
\ref{def_nedc}. Hence, $C_{\alpha}$ has exponential decay.
\end{remark}

\begin{remark}\rm
Sometimes it is useful to consider in Definition \ref{def_neg} or in Definition \ref{def_nedc} the particular case
$\omega(s)\equiv \omega, \ \forall s\geq 0$.
\end{remark}

\section{Exponential stability and instability}

Let $C: T\times Y\rightarrow Y$,
$C(t,s,x,v)=(\varphi(t,s,x),\Phi(t,s,x)v)$ be a skew-evolution semiflow on $Y$.

\begin{definition}\rm\label{def_ns_nes}
The skew-evolution semiflow $C$ is called

$(s)$ \emph{stable} if there exists a mapping
$N:\mathbb{R}_{+}\rightarrow\mathbb{R}_{+}^{*}$ such that
\begin{equation}
\left\Vert \Phi(t,t_{0},x)v\right\Vert \leq N(s)\left\Vert
\Phi(s,t_{0},x)v\right\Vert
\end{equation}%
for all $(t,s),(s,t_{0})\in  T$ and all $(x,v)\in Y$

$(es)$ \emph{exponentially stable} if there exist a mapping
$N:\mathbb{R}_{+}\rightarrow\mathbb{R}_{+}^{*}$ and a constant
$\nu >0$ such that
\begin{equation}
\left\Vert \Phi(t,t_{0},x)v\right\Vert \leq N(s)e^{-\nu
(t-s)}\left\Vert \Phi(s,t_{0},x)v\right\Vert,
\end{equation}%
for all $(t,s),(s,t_{0})\in  T$ and all $(x,v)\in Y$.
\end{definition}

\begin{remark}\rm\label{obs_ns}
The exponential stability of a skew-evolution semiflow implies the stability and, further, the exponential growth.
\end{remark}

In what follows we will present an example of a skew-evolution semiflow that is exponentially
stable but not uniformly exponentially stable.

\begin{example}\rm\label{ex_nues1}
Let $X=\mathbb{R}_{+}$ and $V=\mathbb{R}$. We consider the continuous function
\[
f:\mathbb{R}_{+}\rightarrow[1,\infty), \ f(n)=e^{2n} \
\textrm{and} \ f\left(n+\frac{1}{e^{n^{2}}}\right)=1
\]
and the mapping
\[
\Phi_{f}:
T\times\mathbb{R}_{+}\rightarrow\mathcal{B}(\mathcal{\mathbb{R}}),
\ \Phi_{f}(t,s,x)v=\frac{f(s)}{f(t)}e^{-(t-s)}v.
\]
Then $C_{f}=(\varphi,\Phi_{f})$ is a skew-evolution semiflow on $Y=\mathbb{R}_{+}\times \mathbb{R}$
over all evolution semiflows $\varphi$ on $\mathbb{R}_{+}$. As
\[
\left \Vert \Phi_{f}(t,s,x)v\right\Vert \leq f(s)e^{-(t-s)}|v|, \ \forall
(t,s,x,v)\in T\times Y,
\]
it follows that $C_{f}$ is exponentially stable and, according to Remark \ref{obs_ns}, stable. On the other hand, as
\[
\Phi_{f}\left(n+\frac{1}{e^{n^{2}}},n,x\right)=e^{2n-e^{-n^{2}}}\rightarrow\infty
\ \textrm{when} \ n\rightarrow\infty,
\]
it follows that $C_{f}$ is not uniformly exponentially stable.
\end{example}

\begin{definition}\rm\label{i}\label{cis}
The skew-evolution semiflow $C$ is called

$(is)$ \emph{instable} if there exists a mapping
$N:\mathbb{R}_{+}\rightarrow\mathbb{R}_{+}^{*}$ such that
\begin{equation}
N(t)\left\Vert \Phi(t,t_{0},x)v\right\Vert \geq \left\Vert
\Phi(s,t_{0},x)v\right\Vert
\end{equation}
for all $(t,s),(s,t_{0})\in  T$ and all $(x,v)\in Y$

$(eis)$ \emph{exponentially instable} if there exist a mapping
$N:\mathbb{R}_{+}\rightarrow\mathbb{R}_{+}^{*}$ and a constant
$\nu> 0$ such that
\begin{equation}
N(t)\left\Vert\Phi(t,t_{0},x)v\right\Vert\geq
e^{\nu(t-s)}\left\Vert\Phi(s,t_{0},x)v\right\Vert
\end{equation}
for all $(t,s),(s,t_{0})\in T$ and all $(x,v)\in Y$.
\end{definition}

\begin{remark}\rm\label{obs_nis}
The exponential instability of a skew-evolution semiflow implies
the instability and, also, the exponential decay.
\end{remark}

There exist skew-evolution semiflows that are exponentially
instable but not uniformly exponentially instable, as following example shows.

\begin{example}\rm\label{ex_nues2}
Let $ X=\mathbb{R}_{+}$ and $V=\mathbb{R}$. We consider the
function
\[
f:\mathbb{R}_{+}\rightarrow[1,\infty), \ f(n)=1 \ \textrm{and} \
f\left(n+\frac{1}{e^{n^{2}}}\right)=e^{2n}
\]
and the mapping
\[
\Phi_{f}: T\times\mathbb{R}_{+}\rightarrow\mathcal{B}(
\mathbb{R}), \ \Phi_{f}(t,s,x)v=\frac{f(s)}{f(t)}e^{(t-s)}v.
\]
Then $C_{f}=(\varphi,\Phi_{f})$ is a skew-evolution semiflow on $Y=\mathbb{R}_{+}\times \mathbb{R}$
for all evolution semiflows $\varphi$ on $\mathbb{R}_{+}$. We have
\[
\left \Vert \Phi_{f}(t,s,x)v\right\Vert \geq \frac{1}{f(t)}e^{(t-s)}|v|, \
\forall (t,s,x,v)\in T\times Y,
\]
which proves that $C_{f}$ is exponentially instable and, as in Remark \ref{obs_nis}, instable. But, as
\[
\Phi_{f}\left(n+\frac{1}{e^{n^{2}}},n,x\right)=e^{-2n+e^{-n^{2}}}\rightarrow
0 \ \textrm{for} \ n\rightarrow\infty,
\]
we obtain that $C_{f}$ is not uniformly exponentially instable.
\end{example}

\section{Exponential dichotomy}

Let $C: T\times Y\rightarrow Y$,
$C(t,s,x,v)=(\varphi(t,s,x),\Phi(t,s,x)v)$ be a skew-evolution semiflow on $Y$.

\begin{definition}\rm\label{comp_2}
Two projector families $\{P_{k}\}_{k\in \{1,2\}}$ are said to be
\emph{compatible} with a skew-evolution semiflow
$C=(\varphi,\Phi)$ if

$(dc_{1})$ $P_{1}(x)+P_{2}(x)=I$,
$P_{1}(x)P_{2}(x)=P_{2}(x)P_{1}(x)=0$

$(dc_{2})$
$P_{k}(\varphi(t,s,x))\Phi(t,s,x)v=\Phi(t,s,x)P_{k}(x)v$, $k\in
\{1,2\}$

\noindent for all $t\geq s\geq t_{0}\geq 0$ and all $(x,v)\in Y.$
\end{definition}

\begin{definition}\rm\label{def_d_ed}
The skew-evolution semiflow $C=(\varphi,\Phi)$ is called

$(d)$ \emph{dichotomic} if there exist two projectors $P_{1}$ and
$P_{2}$ compatible with $C$ and some mappings $N_{1}$,
$N_{2}:\mathbb{R}_{+}\rightarrow \mathbb{R}_{+}^{\ast }$ such that
\begin{equation}\label{d_stab}
\left\Vert \Phi(t,t_{0},x)P_{1}(x)v\right\Vert \leq
N_{1}(s)\left\Vert \Phi(s,t_{0},x)P_{1}(x)v\right\Vert
\end{equation}
\begin{equation}\label{d_instab}
\left\Vert \Phi(s,t_{0},x)P_{2}(x)v\right\Vert \leq
N_{2}(t)\left\Vert \Phi(t,t_{0},x)P_{2}(x)v\right\Vert
\end{equation}
for all $(t,s),(s,t_{0})\in T$ and all $(x,v)\in Y$.

$(ed) $\emph{exponentially dichotomic} if there exist two
projectors $P_{1}$ and $P_{2}$ compatible with $C$, some mappings
$N_{1}$, $N_{2}:\mathbb{R}_{+}\rightarrow \mathbb{R}_{+}^{\ast }$
and some constants $\nu_{1}$, $\nu_{2}>0$ such that
\begin{equation}\label{dich_stab}
e^{\nu_{1}(t-s)}\left\Vert \Phi(t,t_{0},x)P_{1}(x)v\right\Vert
\leq N_{1}(s)\left\Vert \Phi(s,t_{0},x)P_{1}(x)v\right\Vert
\end{equation}
\begin{equation}\label{dich_instab}
e^{\nu_{2}(t-s)}\left\Vert \Phi(s,t_{0},x)P_{2}(x)v\right\Vert
\leq N_{2}(t)\left\Vert \Phi(t,t_{0},x)P_{2}(x)v\right\Vert
\end{equation}
for all $(t,s),(s,t_{0})\in T$ and all $(x,v)\in Y$.
\end{definition}

\begin{remark}\rm\label{ed_d}
$(i)$ An exponentially dichotomic skew-evolution semiflow is
dichotomic;

$(ii)$ For $P_{2}=0$ in Definition \ref{def_d_ed} are obtained the
stability, respectively the exponential stability properties for
skew-evolution semiflows;

$(iii)$ For $P_{1}=0$ we obtain in Definition
\ref{def_d_ed} the instability and the exponential instability properties for skew-evolution semiflows.
\end{remark}

\begin{remark}\rm
Without any loss of generality we can consider
\begin{equation*}
N(t)=\max \{N_{1}(t),N_{2}(t)\}, \ t\geq 0 \ \textrm{and} \
\nu=\min \{\nu_{1},\nu_{2}\}.
\end{equation*}
We will call $N_{1}$, $N_{2}$, $\nu_{1}$, $\nu_{2}$, respectively
$N$, $\nu$ the \emph{dichotomic characteristics} asociated to the
skew-evolution semiflow $C$.
\end{remark}

\begin{remark}\rm
Let us consider that the shifted skew-evolution semiflows $C_{\lambda}=(\varphi, \Phi_{\lambda})$ and
$C_{\mu}=(\varphi,\Phi_{\mu})$, where $\Phi_{\lambda}$ and
$\Phi_{\mu}$ are evolution cocycles defined by relation (\ref{relcevshift}) with $\lambda <\mu$,
are exponentially dichotomic with characteristics
\begin{equation*}
N_{\lambda}:\mathbb{R}_{+}\rightarrow \mathbb{R}_{+}^{\ast } \
\textrm{and} \ \nu_{\lambda}>0, \ \textrm{respectively} \
N_{\mu}:\mathbb{R}_{+}\rightarrow \mathbb{R}_{+}^{\ast } \
\textrm{and} \ \nu_{\mu}>0.
\end{equation*}
If we denote
\begin{equation*}
N(t)=\max \{N_{\lambda}(t),N_{\mu}(t)\} \ \textrm{and} \
\nu=\min\{\nu_{\lambda},\nu_{\mu}\}
\end{equation*}
then these are appropriate
for both $C_{\lambda}$ and $C_{\mu}$.
\end{remark}

There exist exponentially dichotomic skew-evolution semiflows that are not
uniformly exponentially dichotomic, as in the next

\begin{example}\rm\label{ex_nued}
Let $ X=\mathbb{R}_{+}$ and $V=\mathbb{R}^{2}$ endowed with the norm
\begin{equation*}
\left\Vert(v_{1},v_{2})\right\Vert=|v_{1}|+|v_{2}|, \
v=(v_{1},v_{2})\in V.
\end{equation*}
The mapping $\Phi: T\times  X\rightarrow \mathcal{B}(V)$,
defined by
\begin{equation*}
\Phi(t,t_{0},x)(v_{1},v_{2})=(e^{t\sin t-s\sin s-2t+2s
}v_{1},e^{2t-2s-3t\cos t+3s\cos s}v_{2})
\end{equation*}
is an evolution cocycle over all evolution semiflows $\varphi$. We
consider the projectors compatible with $C$
\begin{equation*}
P_{1}(x)(v_{1},v_{2})=(v_{1},0) \ \textrm{and} \
P_{2}(x)(v_{1},v_{2})=(0,v_{2}).
\end{equation*}
As
\[
t\sin t-s\sin s-2t+2s\leq -t+3s, \ \forall (t,s)\in T,
\]
we obtain that
\[
\left\Vert\Phi(t,s,x)P_{1}(x)v\right\Vert\leq e^{2s}e^{-(t-s)}|
v_{1}|,\ \forall (t,s,x,v)\in T\times Y.
\]
Similarly, as
\[
2t-2s-3t\cos t+3s\cos s\geq -t-5s, \ \forall (t,s)\in T,
\]
it follows that
\[
e^{6t}\left \Vert\Phi(t,s,x)P_{2}(x)v\right\Vert \geq
e^{5(t-s)}|v_{2}|,\ \forall (t,s,x,v)\in T\times Y.
\]
The skew-evolution semiflow $C=(\varphi,\Phi)$ is exponentially
dichotomic with the dichotomic characteristics
\begin{equation*}
N(t_{0})=e^{6t_{0}} \ \textrm{and} \  \nu=2
\end{equation*}
and, according to Remark \ref{ed_d}, dichotomic. But, as
\[
\Phi\left(2n\pi,2n\pi-\frac{\pi}{2},x\right)=e^{2n\pi-\frac{3\pi}{2}}\rightarrow\infty
\ \textrm{as} \ n\rightarrow\infty,
\]
and
\[
\Phi\left(2n\pi,2n\pi+\frac{\pi}{2},x\right)=e^{-6n\pi-\pi}\rightarrow
0 \ \textrm{as} \ n\rightarrow\infty,
\]
we obtain that $C$ is not uniformly exponentially dichotomic.
\end{example}

In what follows we will denote
\[
C_{k}(t,s,x,v)=(\varphi(t,s,x),\Phi_{k}(t,s,x)v), \ \forall
(t,t_{0},x,v)\in  T\times Y, \ \forall k\in \{1,2\},
\]
where
\[
\Phi_{k}(t,t_{0},x)=\Phi(t,t_{0},x)P_{k}(x), \ \forall
(t,t_{0})\in  T, \ \forall x\in  X, \ \forall k\in \{1,2\}.
\]

We give an integral characterization for the dichotomy property of skew-evolution semiflows.

\begin{theorem}\label{tned}
A strongly measurable skew-evolution semiflow $C=(\varphi,\Phi)$ is exponentially
dichotomic if and only if there exist two projectors $P_{1}$ and $P_{2}$ compatible with $C$
such that $C_{1}$ has exponential growth and $C_{2}$ has exponential decay,
some functions $M_{1}$, $M_{2}:\mathbb{R}_{+}\rightarrow\mathbb{R}_{+}^{*}$ and
some constants $\alpha$, $\beta>0$ such that

$(ed_{1})$
\begin{equation}\label{dich_stab_int}
\int_{t_{0}}^{t}e^{\alpha (\tau-t_{0})}\left\Vert \Phi_{1}(\tau
,t_{0},x)v\right\Vert d\tau \leq M_{1}(t_{0})\left\Vert
P_{1}(x)v\right\Vert
\end{equation}

$(ed_{2})$
\begin{equation}\label{dich_instab_int}
\int\limits_{t_{0}}^{t}e^{\beta(t-\tau)}\left\Vert\Phi_{2}(\tau,t_{0},x)v\right\Vert
d\tau \leq M_{2}(t)\left\Vert\Phi_{2}(t,t_{0},x)v\right\Vert
\end{equation}
for all $(t,t_{0})\in T$ and all $(x,v)\in Y$.
\end{theorem}

\begin{proof}
\textit{Necessity}. As the skew-evolution semiflow $C$ is
exponentially dichotomic, there exist two projectors $P_{1}$ and
$P_{2}$ compatible with $C$, a function
$N_{1}:\mathbb{R}_{+}\rightarrow \mathbb{R}_{+}^{\ast }$ and a
constant $\nu_{1} >0$ such that
\begin{equation*}
\left\Vert \Phi(t,s,x)P_{1}(x)v\right\Vert \leq
N_{1}(s)e^{-\nu_{1} (t-s)}\left\Vert P_{1}(x)v\right\Vert
\end{equation*}%
for all $t\geq s\geq 0$ and all $(x,v)\in Y$. We consider $\alpha$
such that $\nu_{1}=2\alpha$. Following relations hold
\begin{equation*}
\int_{t_{0}}^{t}e^{\alpha (\tau-t_{0})}\left\Vert
\Phi(\tau,t_{0},x)v\right\Vert d\tau \leq N_{1}(t_{0})\left\Vert
v\right\Vert \int_{t_{0}}^{t}e^{\alpha (\tau-t_{0})}e^{-\nu_{1}
(\tau -t_{0})} d\tau \leq M_{1}(t_{0})\left\Vert v\right\Vert,
\end{equation*}%
where we have denoted
\begin{equation*}
M_{1}(t_{0})=\frac{N_{1}(t_{0})}{\alpha}.
\end{equation*}
Hence, relation (\ref{dich_stab_int}) is obtained.

Also, there exist a function
$N_{2}:\mathbb{R}_{+}\rightarrow\mathbb{R}_{+}^{*}$ and a constant
$\nu_{2}> 0 $ such that
\begin{equation*}
e^{\nu_{2}(t-s)}\left\Vert\Phi(s,t_{0},x)P_{2}(x)v\right\Vert \leq
N_{2}(t)\left\Vert\Phi(t,t_{0},x)P_{2}(x)v\right\Vert
\end{equation*}
for all $t\geq t_{0}\geq 0$ and all $(x,v)\in Y$. We will consider
$\beta$ such that $\nu_{2}=2\beta$. We have
\begin{equation*}
\int_{t_{0}}^{t}e^{\beta(t-\tau) }\left\Vert \Phi(\tau
,t_{0},x)P_{2}(x)v\right\Vert d\tau \leq
\end{equation*}
\begin{equation*}
\leq N_{2}(t)\int_{t_{0}}^{t}e^{\beta(t-\tau) }e^{-\nu_{2}(t-\tau
)}\left\Vert \Phi(t,t_{0},x)P_{2}(x)v\right\Vert d\tau\leq
M_{2}(t)\left\Vert \Phi(t,t_{0},x)P_{2}(x)v\right\Vert,
\end{equation*}
where we have denoted
\begin{equation*}
M_{2}(t)=\frac{N_{2}(t)}{\beta}.
\end{equation*}
Relation (\ref{dich_instab_int}) is then obtained.

\textit{Sufficiency.} As $C_{1}$ has exponential growth, similarly
as in Theorem 2.3 of \cite{MeStBu_UVT}, proved for evolution
operators, there exists a nondecreasing function $f:[0,\infty
)\rightarrow [1,\infty )$ with the property
$\underset{t\rightarrow\infty}\lim f(t)=\infty$ such that
\begin{equation*}
\left\Vert \Phi(t,t_{0},x)v\right\Vert \leq f(t-s)\left\Vert
\Phi(s,t_{0},x)v\right\Vert,
\end{equation*}%
for all $(t,s),(s,t_{0})\in  T$ and all $(x,v)\in Y$.

For $t\geq t_{0}+1$ we have
\[
\left\Vert \Phi(t,t_{0},x)P_{1}(x)v\right\Vert e^{\alpha (t-t_{0})} \int_{0}^{1}\frac{e^{-\alpha u}}{f(u)}%
du=
\]
\[
=e^{\alpha (t-t_{0})}\left\Vert
\Phi(t,t_{0},x)P_{1}(x)v\right\Vert\int_{t-1}^{t}\frac{e^{-\alpha
(t-\tau)}}{f(t-\tau)}d\tau=
\]
\[
=\int_{t-1}^{t}\frac{\left\Vert\Phi(t,t_{0},x)P_{1}(x)v\right\Vert}{f(t-\tau)}e^{\alpha
(\tau-t_{0})}d\tau
\leq\int_{t-1}^{t}\left\Vert\Phi(\tau,t_{0},x)P_{1}(x)v\right\Vert
e^{\alpha (\tau-t_{0})}d\tau\leq
\]
\[
\leq M_{1}(t_{0})\left\Vert P_{1}(x)v\right\Vert.
\]
For $t\in [t_{0},t_{0}+1)$ we have
\begin{equation*}
\left\Vert \Phi_{-\alpha}(t,t_{0},x)P_{1}(x)v\right\Vert \leq
f(1)e^{\alpha }\left\Vert P_{1}(x)v\right\Vert,
\end{equation*}%
where the evolution cocycle $\Phi_{-\alpha}$ is given by relation
(\ref{relcevshift}). We obtain
\begin{equation*}
\left\Vert \Phi_{-\alpha }(t,t_{0},x)P_{1}(x)v\right\Vert \leq N_{1}%
(t_{0})\left\Vert P_{1}(x)v\right\Vert,
\end{equation*}%
for all $(t,t_{0})\in T$ and all $(x,v)\in Y$, where we have denoted
\begin{equation*}
N_{1}(t_{0})=f(1)e^{\alpha}+M_{1}(t_{0})\left[\int_{0}^{1}\frac{e^{-\alpha
u}}{f(u)}du\right]^{-1}.
\end{equation*}
It follows that
\begin{equation*}
\left\Vert \Phi(t,t_{0},x)P_{1}(x)v\right\Vert \leq N_{1}%
(t_{0})e^{-\alpha (t-t_{0})}\left\Vert P_{1}(x)v\right\Vert,
\end{equation*}%
for all $t\geq t_{0}\geq 0$ and all $(x,v)\in Y$. Hence, relation (\ref{dich_stab}) was proved.

As $C_{2}$ has exponential decay, by a similar deduction used to
prove Theorem 3.3 of \cite{MeStBu_UVT} for evolution operators,
there exists a nondecreasing function $g:[0,\infty )\rightarrow
[1,\infty )$ with the property $\underset{t\rightarrow\infty}\lim
g(t)=\infty$ such that
\begin{equation*}
\left\Vert \Phi(s,t_{0},x)v\right\Vert \leq g(t-s)\left\Vert
\Phi(t,t_{0},x)v\right\Vert,
\end{equation*}%
for all $(t,s),(s,t_{0})\in  T$ and all $(x,v)\in Y$.

For $t\geq s\geq t_{0}\geq 0$, $x\in  X$, $v\in V$ we have
\[
\left\Vert \Phi(s,t_{0},x)P_{2}(x)v\right\Vert e^{\beta (t-s)} \int_{0}^{1}\frac{e^{\beta u}}{g(u)}%
du=
\]
\[
=\int_{s-1}^{s}\left\Vert
\Phi(s,t_{0},x)P_{2}(x)v\right\Vert\frac{e^{\beta
(t-s)}}{g(s-\tau)}e^{\beta (s-\tau)}d\tau\leq
\]
\[
\leq\int_{s-1}^{s}\left\Vert\Phi(\tau,t_{0},x)P_{2}(x)v\right\Vert
e^{\beta (t-\tau)}d\tau\leq M_{2}(t)\left\Vert
\Phi(t,t_{0},x)P_{2}(x)v\right\Vert.
\]
We obtain
\begin{equation*}
\left\Vert \Phi(s,t_{0},x)P_{2}(x)v\right\Vert \leq
N_{2}(t)e^{-\beta (t-s)}\left\Vert
\Phi(t,t_{0},x)P_{2}(x)v\right\Vert,
\end{equation*}
for all $t\geq t_{0}\geq 0$ and all $(x,v)\in Y$, where we have denoted
\begin{equation*}
N_{2}(t)=M_{2}(t) \left[\int_{0}^{1}\frac{e^{\beta
u}}{g(u)}\right]^{-1}.
\end{equation*}
Relation (\ref{dich_instab}) is then obtained.

Hence, the skew-evolution semiflow $C$ is exponentially dichotomic.
\end{proof}

\section{Exponential trichotomy}

Let $C: T\times Y\rightarrow Y$,
$C(t,s,x,v)=(\varphi(t,s,x),\Phi(t,s,x)v)$ be a skew-evolution semiflow on $Y$.

\begin{definition}\rm\label{comp_3}
Three projector families $\{P_{k}\}_{k\in \{1,2,3\}}$ are said to
be \emph{compatible} with a skew-evolution semiflow
$C=(\varphi,\Phi)$ if

$(tc_{1})$ $P_{1}(x)+P_{2}(x)+P_{3}(x)=I$,
$P_{i}(x)P_{j}(x)=P_{j}(x)P_{i}(x)=0$, $\forall i,j\in \{1,2,3\}$,
$i\neq j$

$(tc_{2})$
$P_{k}(\varphi(t,s,x))\Phi(t,s,x)v=\Phi(t,s,x)P_{k}(x)v$, $\forall
k\in \{1,2,3\}$,

\noindent for all $(t,s),(s,t_{0})\in T$ and all $(x,v)\in Y.$
\end{definition}

\begin{definition}\rm\label{net}
A skew-evolution semiflow $C=(\varphi,\Phi)$ is said to be

$(t)$ \emph{trichotomic} if there exist three projectors $P_{1}$,
$P_{2}$ and $P_{3}$ compatible with $C$ and some functions
$N_{1}$, $N_{2}$, $N_{3}:\mathbb{R}_{+}\rightarrow
\mathbb{R}_{+}^{\ast }$ such that
\begin{equation}
\left\Vert \Phi(t,t_{0},x)P_{1}(x)v\right\Vert \leq
N_{1}(t_{0})\left\Vert \Phi(s,t_{0},x)P_{1}(x)v\right\Vert
\end{equation}
\begin{equation}
\left\Vert \Phi(s,t_{0},x)P_{2}(x)v\right\Vert \leq
N_{2}(t_{0})\left\Vert \Phi(t,t_{0},x)P_{2}(x)v\right\Vert
\end{equation}
\[
\left\Vert \Phi(s,t_{0},x)P_{3}(x)v\right\Vert \leq
N_{3}(t_{0})\left\Vert \Phi(t,t_{0},x)P_{3}(x)v\right\Vert\leq
\]
\begin{equation}
\leq N_{3}^{2}(t_{0})\left\Vert
\Phi(s,t_{0},x)P_{3}(x)v\right\Vert
\end{equation}
for all $(t,s),(s,t_{0})\in T$ and all $(x,v)\in Y$.

$(et)$ \emph{exponentially trichotomic} if there exist three
projectors $P_{1}$, $P_{2}$ and $P_{3}$ compatible with $C$, some
functions $N_{1}$, $N_{2}$, $N_{3}$,
$N_{4}:\mathbb{R}_{+}\rightarrow \mathbb{R}_{+}^{\ast }$ and some
constants $\nu_{1}$, $\nu_{2}$, $\nu_{3}$, $\nu_{4}$ with the
properties $\nu_{1}<\nu_{2}\leq 0\leq \nu_{3}<\nu_{4}$ such that
\begin{equation}\label{Pes}
\left\Vert \Phi(t,t_{0},x)P_{1}(x)v\right\Vert \leq
N_{1}(s)\left\Vert \Phi(s,t_{0},x)P_{1}(x)v\right\Vert
e^{\nu_{1}(t-s)}
\end{equation}
\begin{equation}\label{Qeis}
N_{4}(t)\left\Vert \Phi(t,t_{0},x)P_{2}(x)v\right\Vert \geq
\left\Vert \Phi(s,t_{0},x)P_{2}(x)v\right\Vert e^{\nu_{4}(t-s)}
\end{equation}
\begin{equation}\label{Reg}
\left\Vert \Phi(t,t_{0},x)P_{3}(x)v\right\Vert \leq
N_{3}(s)\left\Vert \Phi(s,t_{0},x)P_{3}(x)v\right\Vert
e^{\nu_{3}(t-s)}
\end{equation}
\begin{equation}\label{Redc}
N_{2}(t)\left\Vert \Phi(t,t_{0},x)P_{3}(x)v\right\Vert \geq
\left\Vert \Phi(s,t_{0},x)P_{3}(x)v\right\Vert e^{\nu_{2}(t-s)}
\end{equation}
for all $(t,s),(s,t_{0})\in T$ and all $(x,v)\in Y$.
\end{definition}

\begin{remark}\rm\label{et_t}
$(i)$ An exponentially trichotomic skew-evolution semiflow is
trichotomic;

$(ii)$ For $P_{1}=0$ in Definition \ref{net} are obtained the
properties of dichotomy, respectively exponential dichotomy;

$(iii)$ For $P_{2}=P_{3}=0$ the properties of stability,
respectively exponential stability follow from Definition
\ref{net};

$(iv)$ If $P_{1}=P_{3}=0$ the properties of instability and exponential instability are obtained from Definition
\ref{net}.
\end{remark}

\begin{remark}\rm
Without any loss of generality we can chose
\begin{equation*}
N(t)=\max\{N_{1}(t), \ N_{2}(t), \ N_{3}(t), \ N_{4}(t)\}, \ t\geq
0
\end{equation*}
respectively
\begin{equation*}
\alpha=-\nu_{1}=\nu_{4}>0 \ \textrm{and} \
\beta=-\nu_{2}=\nu_{3}>0.
\end{equation*}
We call $N_{1}$, $N_{2}$, $N_{3}$, $N_{4}$, $\nu_{1}$, $\nu_{2}$,
$\nu_{3}$, $\nu_{4}$, respectively $N$, $\alpha$, $\beta$ the trichotomic characteristics of the skew-evolution semiflow $C$.
\end{remark}

\begin{example}\rm
Let us consider $( X,d)$ a metric space given as in Example \ref{ex_ses}. The mapping
\begin{equation*}
\varphi: T\times  X\rightarrow  X, \
\varphi(t,s,x)(\tau)=x(t-s+\tau)
\end{equation*}
is an evolution semiflow on $ X$.

Let $V=\mathbb{R}^{3}$ with the norm
\begin{equation*}
\left\Vert(v_{1},v_{2},v_{3})\right\Vert=|v_{1}|+|v_{2}|+|v_{3}|.
\end{equation*}
The mapping $\Phi:
T\times  X\rightarrow \mathcal{B}(V)$, given by
\begin{equation*}
\Phi(t,s,x)(v)=\left(e^{-2(t-s)f(0)+\int_{0}^{t}x(\tau)d\tau}v_{1},
e^{t-s+\int_{0}^{t}x(\tau)d\tau}v_{2},e^{-(t-s)f(0)+2\int_{0}^{t}x(\tau)d\tau}v_{3}\right)
\end{equation*}
is an evolution cocycle.

We consider the projections
\begin{equation*}
P_{1}(x)(v)=(v_{1},0,0), \ P_{2}(x)(v)=(0,v_{2},0), \
P_{3}(x)(v)=(0,0,v_{3}).
\end{equation*}
The skew-evolution semiflow $C=(\varphi,\Phi)$ is exponentially trichotomic
with asociated trichotomic characteristics
\begin{equation*}
\nu_{1}=\nu_{2}=-f(0),  \ \nu_{3}=f(0) \ \textrm{and} \
\nu_{4}=1
\end{equation*}
\begin{equation*}
N_{1}(t_{0})=e^{t_{0}f(0)}, \ \ N_{2}(t)=e^{-2lt}, \ \,
N_{3}(t_{0})=e^{2t_{0}f(0)} \ \textrm{and} \ N_{4}(t)=e^{-lt}.
\end{equation*}
\end{example}

A characterization for the property of exponential trichotomy can
be given by the next

\begin{theorem}\label{tnet}
Let $C=(\varphi,\Phi)$ be a skew-evolution semiflow and three
projectors $P_{1}$, $P_{2}$ and $P_{3}$ compatible with $C$ such
that $C_{1}$ has exponential growth and is $*$-strongly measurable
and $C_{2}$ has exponential decay and is strongly measurable. Then
$C$ is exponentially trichotomic if and only if there exist a
mapping
$\widetilde{N}:\mathbb{R}_{+}\rightarrow\mathbb{R}_{+}^{*}$ and a
constant $\alpha>0$, a mapping
$\overline{N}:\mathbb{R}_{+}\rightarrow\mathbb{R}_{+}^{*}$ and a
constant $\beta>0$, some functions $\widetilde{M}$,
$\overline{M}:[0,\infty )\rightarrow (0,\infty )$, some
nondecreasing functions $\widetilde{g}$, $\overline{g}:[0,\infty
)\rightarrow (0,\infty )$ with the property
$\underset{t\rightarrow\infty}\lim
\widetilde{g}(t)=\underset{t\rightarrow\infty}\lim
\overline{g}(t)=\infty$ such that

$(et_{1})$
\begin{equation}
\int^{t}_{t_{0}}e^{\alpha(t-s)}\left\Vert\Phi(t,s,\varphi(s,t_{0},x))^{*}P_{1}(x)v^{*}\right\Vert
ds\leq \widetilde{N}(t_{0})\left\Vert P_{1}(x)v^{*}\right\Vert,
\end{equation}
for all $(t,t_{0})\in T$ and all $(x,v^{*})\in X\times
V^{*}$ with $\left\Vert v^{*}\right\Vert\leq 1$

$(et_{2})$
\begin{equation}
\int^{t}_{0} e^{-\beta
(s-t_{0})}\left\Vert\Phi(s,t_{0},x)P_{2}(x)v\right\Vert ds\leq
\overline{N}(t)e^{-\beta(t-t_{0})}\left\Vert
\Phi(t,t_{0},x)P_{2}(x)v\right\Vert,
\end{equation}
for all $(t,t_{0})\in T$ and all $(x,v)\in Y$

$(et_{3})$
\begin{equation}
\left\Vert \Phi(t,t_{0},x)P_{3}(x)v\right\Vert \leq
\widetilde{M}(s)\widetilde{g}(t-s)\left\Vert \Phi(s,t_{0},x)P_{3}(x)v\right\Vert
\end{equation}%
for all $(t,s),(s,t_{0})\in T$ and all $(x,v)\in Y$

$(et_{4})$
\begin{equation}
\left\Vert \Phi(s,t_{0},x)P_{3}(x)v\right\Vert \leq
\overline{M}(t)\overline{g}(t-s)\left\Vert \Phi(t,t_{0},x)P_{3}(x)v\right\Vert
\end{equation}%
for all $(t,s,),(s,t_{0})\in T$ and all $(x,v)\in Y.$
\end{theorem}

\begin{proof}
\emph{Necessity.} $(et_{1})$ We consider
\[
\alpha=-\frac{\nu_{1}}{2}
\]
and we obtain
\[
\int^{t}_{t_{0}}e^{\alpha(t-s)}\left\Vert\Phi(t,s,x)^{*}P_{1}(x)v^{*}\right\Vert
ds\leq\int^{t}_{t_{0}}N_{1}(t_{0})e^{\frac{\nu_{1}}{2}(t-s)}\left\Vert
P_{1}(x)v^{*}\right\Vert ds\leq
\]
\[
\leq \widetilde{N}(t_{0})\left\Vert
P_{1}(x)v^{*}\right\Vert,
\]
for all $(t,t_{0},x,v^{*})\in T\times X\times V^{*}$, where we have denoted
\[
\widetilde{N}(t_{0})=-\frac{2}{\nu}N_{1}(t_{0}).
\]

$(et_{2})$ Let us define
\[
\beta=\frac{\nu}{2}>0,
\]
where the existence of constant $\nu$ is assured by hypothesis and
by Definition \ref{cis}. Hence, we obtain
\begin{equation*}
\int_{0}^{t}e^{-\beta(s-t_{0})}\left\Vert
\Phi(s,t_{0},x)v\right\Vert ds \leq
\end{equation*}
\begin{equation*}
\leq N(t)\int_{t_{0}}^{t}e^{-\beta(s-t_{0})}\left\Vert
\Phi(t,t_{0},x)v\right\Vert e^{-\nu (t-s)}ds\leq
\beta^{-1}N(t)e^{-\beta(t-t_{0})}\left\Vert
\Phi(t,t_{0},x)v\right\Vert.
\end{equation*}

$(et_{3})$ It is obtained immediately if we consider
\[
\widetilde{M}(u)=N_{3}(u),\ u\geq 0 \ \textrm{and} \ \widetilde{g}(v)=e^{\nu_{3}
v}, \ v\geq 0.
\]

$(et_{4})$ It follows for
\[
\overline{M}(u)=N_{2}(u),\ u\geq 0 \ \textrm{and} \ \overline{g}(v)=e^{-\nu_{2}
v}, \ v\geq 0.
\]

\emph{Sufficiency}. $(et_{1})$ Let $t\geq t_{0}+1$ and
$s\in[t_{0},t_{0}+1)$. Then
\[
e^{-[\alpha+\omega(t_{0})]}\left|\left\langle
v^{*},e^{\alpha(t-t_{0})}\Phi(t,t_{0},x)v\right\rangle \right|=
\]
\[
=
e^{-[\alpha+\omega(t_{0})]}\int_{t_{0}}^{t_{0}+1}\left|\left\langle
\Phi(t,s,x)^{*}v^{*},e^{\alpha(t-t_{0})}\Phi(s,t_{0},x)v\right\rangle
\right|ds\leq
\]
\[
\leq\int_{t_{0}}^{t_{0}+1}e^{\alpha(t-s)}\left\Vert\Phi(t,s,\varphi(s,t_{0},x))^{*}v^{*}\right\Vert
e^{-\omega(t_{0})(s-t_{0})}\left\Vert\Phi(s,t_{0},x)v\right\Vert
ds\leq
\]
\[
\leq M(t_{0})\left\Vert
v\right\Vert\int_{t_{0}}^{t}e^{\alpha(t-s)}\left\Vert\Phi(t,s,\varphi(s,t_{0},x))^{*}v^{*}\right\Vert
ds\leq M(t_{0})N(t_{0})\left\Vert v\right\Vert\left\Vert
v^{*}\right\Vert,
\]
where the functions $M$,
$\omega:\mathbb{R}_{+}\rightarrow\mathbb{R}_{+}^{*}$ are given by
Definition \ref{def_neg}. By taking supremum over $\left\Vert
v^{*}\right\Vert\leq 1$ we obtain
\[
e^{-[\alpha+\omega(t_{0})]}e^{\alpha(t-t_{0})}\left\Vert\Phi(t,t_{0},x)v\right\Vert\leq
M(t_{0})N(t_{0})\left\Vert v\right\Vert, \ \forall t\geq t_{0}+1,
\]
and, further,
\[
\left\Vert\Phi(t,t_{0},x)v\right\Vert\leq
M_{1}(t_{0})e^{-\alpha(t-t_{0})}\left\Vert v\right\Vert,\ \forall
t\geq t_{0}+1,
\]
where we have denoted
\[
M_{1}(t_{0})=M(t_{0})N(t_{0})e^{[\alpha+\omega(t_{0})]}, \
t_{0}\geq 0.
\]
For $t\in[t_{0},t_{0}+1)$ we have
\[
\left\Vert\Phi(t,t_{0},x)v\right\Vert\leq M(t_{0})e^{\omega
(t_{0})(t-t_{0})}\left\Vert v\right\Vert\leq
M_{2}(t_{0})e^{-\alpha(t-t_{0})}\left\Vert v\right\Vert,
\]
where we have denoted
\[
M_{2}(t_{0})=M(t_{0})e^{[\alpha+\omega(t_{0})]}, \ t_{0}\geq 0.
\]
Hence
\[
\left\Vert\Phi(t,t_{0},x)v\right\Vert\leq[M_{1}(t_{0})+M_{2}(t_{0})]e^{-\alpha(t-t_{0})}\left\Vert
v\right\Vert,\ \forall (t,t_{0},x,v)\in T\times X\times V,
\]
which proves relation (\ref{Pes}) of Definition \ref{et_t}.

$(et_{2})$ We denote $K=\int\limits_{0}^{1}e^{-\beta%
u}f(u)du$, where function $f$ is given as in Theorem 3.3 of
\cite{MeStBu_UVT}. We obtain succesively
\begin{equation*}
K\left\Vert v\right\Vert=\int_{t_{0}}^{t_{0}+1}e^{-\beta
(\tau-t_{0})} f(\tau -t_{0})\left\Vert
\Phi(t_{0},t_{0},x)v\right\Vert d\tau \leq
\end{equation*}%
\begin{equation*}
\leq \int_{t_{0}}^{t_{0}+1}e^{-\beta (\tau-t_{0})}\left\Vert
\Phi(\tau,t_{0},x)v\right\Vert d\tau \leq M(t)\left\Vert
\Phi_{\beta}(t,t_{0},x)v\right\Vert=
\end{equation*}
\[
=M(t)e^{-\beta (t-t_{0})}\left\Vert \Phi (t,t_{0},x)v\right\Vert
\]
for all $(t,t_{0})\in  T$ and all $(x,v)\in Y$. Hence, relation
(\ref{Qeis}) of Definition \ref{et_t} was proved.

$(et_{3})$ As $\underset{t\rightarrow\infty}\lim
\widetilde{g}(t)=\infty$, there exists $\delta>0$ such that $\widetilde{g}(\delta)>1$.
Let $(t,s)\in T$. Then there exists $n\in\mathbb{N}$ and
$r\in[0,\delta)$ such that $t=s+n\delta+r$. We obtain succesively
\[
\left\Vert \Phi(t,t_{0},x)P_{3}(x)v\right\Vert=\left\Vert
\Phi(s+n\delta+r,t_{0},x)P_{3}(x)v\right\Vert\leq
\]
\[
\leq \widetilde{g}(r) \widetilde{M}(s+n\delta)\left\Vert
\Phi(s+n\delta,t_{0},x)P_{3}(x)v\right\Vert\leq
\]
\[
\leq \widetilde{g}(\delta)
\widetilde{M}(s+n\delta)\left\Vert
\Phi(s+n\delta,t_{0},x)P_{3}(x)v\right\Vert\leq
\]
\[
\leq [\widetilde{g}(\delta)]^{2}\widetilde{M}(s+(n-1)\delta+r)\left\Vert
\Phi(s+(n-1)\delta,t_{0},x)P_{3}(x)v\right\Vert\leq...\leq
\]
\[
\leq
[\widetilde{g}(\delta)]^{n+1}\widetilde{M}(s+(n-1)\delta)...\widetilde{M}(s)\left\Vert
\Phi(s+r,t_{0},x)P_{3}(x)v\right\Vert,
\]
for all $(t,s),\ (s,t_{0})\in T$ and all $(x,v)\in Y$.
If we define
\[
N_{3}(u)=\widetilde{g}(\delta)\widetilde{M}(u+(n-1)\delta+r)...\widetilde{M}(u+r)\widetilde{M}(u)
\ \textrm{and} \ \nu_{3}=\frac{\ln \widetilde{g}(\delta)}{\delta}
\]
we obtain relation (\ref{Reg}) of Definition \ref{et_t}.

$(et_{4})$ Without loss of generality we can consider $\overline{g}(1)>1$. Let $(t,t_{0})\in T$ and
$n=[t-t_{0}]$. We obtain $t_{0}+n\leq t< t_{0}+n+1$. Following relations hold for all $(x,v)\in
Y$
\[
\left\Vert \Phi(t,t_{0},x)P_{3}(x)v\right\Vert\geq
\frac{1}{\overline{M}(t)}\frac{1}{\overline{g}(1)}\left\Vert
\Phi(t-1,t_{0},x)P_{3}(x)v\right\Vert\geq
\]
\[
\geq\frac{1}{\overline{M}(t)}\frac{1}{\overline{M}(t-1)}\frac{1}{[\overline{g}(1)]^{2}}\left\Vert
\Phi(t-2,t_{0},x)P_{3}(x)v\right\Vert\geq...\geq
\]
\[
\geq\frac{1}{\overline{M}(t)}\frac{1}{\overline{M}(t-1)}...\frac{1}{\overline{M}(t-(n-1))}\frac{1}{[\overline{g}(1)]^{n}}\left\Vert
\Phi(t-n,t_{0},x)P_{3}(x)v\right\Vert\geq
\]
\[
\geq\frac{1}{\overline{M}(t)}\frac{1}{\overline{M}(t-1)}...
\frac{1}{\overline{M}(t-(n-1))}\frac{1}{\overline{M}(t-t_{0}-n)}\frac{1}
{[\overline{g}(1)]^{n}}\frac{1}{\overline{g}(t-t_{0}-n)}\left\Vert P_{3}(x)
v\right\Vert.
\]
If we denote
\[
N_{2}(t)=\overline{g}(1)\overline{M}(t)\overline{M}(t-1)...\overline{M}(1)
\ \textrm{and} \ \nu_{2} =\overline{g}(1),
\]
relation (\ref{Redc}) of Definition \ref{et_t} is obtained.

Hence, the skew-evolution semiflow $C$ is exponentially trichotomic.
\end{proof}

\footnotesize{

\vspace{5mm}

\noindent\begin{tabular}[t]{ll}

\textsc{Codru\c{t}a Stoica} \\
Institut de Math\' ematiques  \\
Universit\' e Bordeaux 1  \\
351 Cours de la Libération \\
F-33405  Talence Cedex
France  \\
e-mail: \texttt{codruta.stoica@math.u-bordeaux1.fr}

\vspace{3mm}\\

\textsc{Mihail Megan}\\
Faculty of Mathematics\\
West University of Timi\c{s}oara\\
Bd. P\^{a}rvan, No.4\\
300223 Timi\c{s}oara\\
Romania\\
E-mail: \texttt{megan@math.uvt.ro}
\end{tabular}
}

\end{document}